\documentclass[12pt,a4paper]{amsart}

\numberwithin{equation}{subsection}

\newtheorem{theorem}[subsection]{Theorem}
\newtheorem{lemma}[subsection]{Lemma}
\newtheorem{proposition}[subsection]{Proposition}
\newtheorem{corollary}[subsection]{Corollary}

\newtheorem{Sketch of the proof}[subsection]{Sketch of the proof}
\theoremstyle{definition}

\theoremstyle{remark}

\newcommand{\zp}{{{\bf{Z}}/p}}

\newcommand{\fq}{{\bf{F}}_q}
\newcommand{\fp}{{\bf{F}}_p}
\newcommand{\field}{{\bf{F}}}

\newcommand{\fieldv}{{\field}[V]}

\newcommand{\tr}{\mathop{\rm Tr}}
\newcommand{\Tr}{\mathop{\rm Tr}}

\newcommand{\lt}{\mathop{\rm LT}}

\title[Noether numbers]{The Noether numbers for cyclic groups\\ of prime order}

\author{P. Fleischmann}
\author{M. Sezer}
\author{R.J. Shank}
\author{C.F. Woodcock}

\address{Department of Mathematics, Bo\u gazi\c ci University \\
\hfil\break\indent TR-34342 Bebek, Istanbul, Turkey}

\email{mufit.sezer@boun.edu.tr}

\address{Institute of Mathematics, Statistics \&  Actuarial Science \\
 \hfil\break\indent University of Kent at Canterbury, CT2 7NF, UK}

\email{P.Fleischmann@kent.ac.uk}
\email{R.J.Shank@kent.ac.uk}
\email{C.F.Woodcock@kent.ac.uk}

\thanks{Research supported by EPSRC}

\date{\today}

\begin{document}

\begin{abstract}
The Noether number of a representation is the largest degree of an element 
in a minimal homogeneous generating set for the corresponding ring of invariants.
We compute the Noether number for an arbitrary representation of a cyclic group of prime order,
and as a consequence prove the ``$2p-3$ conjecture''.
\end{abstract}

\maketitle

\section{Introduction}

Let $V$ denote a finite dimensional representation of a finite group $G$ over a field $\field$, 
i.e., $V$ is a finite  module over the group ring $\field G$.
The action of $G$ on $V$ induces an action of $G$ on the dual $V^*$ which extends to an action by
degree preserving algebra automorphisms on the symmetric algebra of $V^*$, $\field[V]:=S(V^*)$.
The {\it ring of invariants},
$$\field[V]^G:=\{f\in \field[V]\mid g(f)=f,\ \forall g\in G\},$$
is a finitely generated subring of $\field[V]$, and can be interpreted as the ring
of regular functions on the categorical quotient $V/\!/G$. Since the action of $G$ on $\field[V]$ preserves degree,
the generators can be chosen to be homogeneous. The {\it Noether number} of the representation, $\beta(V)$,
is defined to be the least integer $d$ such that $\fieldv^G$ is generated by homogeneous elements of degree 
less than or equal to $d$.
Emmy Noether \cite{Noether} proved that if $\field$ has characteristic zero,
then $\beta(V)\leq |G|$. This result has been extended to all non-modular representations by 
Fleischmann \cite{Fleischmann} and Fogarty \cite{Fogarty}. 
It is clear that a non-trivial non-modular representation of a cyclic group of prime order has
Noether number $|G|$. In fact, cyclic groups are the only groups with non-modular representations having
$|G|$ as the Noether number (see Schmid \cite{schmid} for characteristic zero and Sezer \cite{sezer} for the
generalisation to non-modular representations). 
It follows from the work of Richman (see \cite{richman1} and \cite{richman2}),
that $|G|$ is not an upper bound for $\beta(V)$ for modular representations.
It has been conjectured that $\beta(V)\leq {\rm max}\{|G|,\dim(V)(|G|-1)\}$ (see \cite[Remark~3.9.10]{de-ke}).
Derksen \& Kemper have shown that
$$\beta(V)\leq n(|G|-1)+|G|^{n\cdot 2^{n-1}+1}\cdot n^{2^{n-1}+1}$$ 
with $n:=\dim(V)$ \cite[Theorem~3.9.11]{de-ke}, and
Karagueuzian \& Symonds  \cite{ka-sy} have shown that if $\field=\fq$ then
$\beta(V)\leq (q^n-1)(nq-n-1)/(q-1)$ (for $n>2$).
The Noether number is an important characteristic of the ring of invariants.
For example, {\it a prior} knowledge of the Noether number reduces the construction of a generating set
to a problem in linear algebra.

For the remainder of the paper, $\field$ will denote a field of characteristic $p$ for a prime number $p$,
and $\zp$ will denote the cyclic group of order $p$. Despite the simplicity of the representation theory of $\zp$, 
the corresponding rings of invariants can be surprisingly complicated.
From Ellingsrud \& Skjelbred \cite{el-skj}, we know that these rings of invariants almost always fail to be Cohen-Macaulay.
While formulae for the Hilbert series are known for all representations of $\zp$, see Almkvist \& Fossum \cite{alm-fos} and
Hughes \& Kemper \cite{hk},
explicit generating invariants are known for only a few special cases. 

Choose a generator $\sigma$ for $\zp$.
In $\field\zp$, define $\Delta:=\sigma-1$ and $\Tr:=\sum_{i=1}^p\sigma^i$.
The isomorphism
type of a representation of $\zp$ is determined by the Jordan canonical form of $\sigma$. 
If $n\leq p$, then the $n\times n$ matrix over $\field$ consisting of a single Jordan block with eigenvalue
$1$, has order $p$ and determines an indecomposable $\field\zp$-module, which we denote by $V_n$.
Note that for $n>p$, this matrix has order at least $p^2$ and, therefore, does not give a representation of $\zp$.
Recall that $1$ is the only $p^{th}$ root of unity in $\field$. Thus $V_1, V_2,\ldots,V_p$ is a complete list, up to isomorphism,
of the indecomposable $\field\zp$-modules. 
For a positive integer $k$, we will denote the 
direct sum of  $k$ copies of $V_n$ by $kV_n$. 
Explicit generating sets for $\field[V_2]^{\zp}$ and $\field[V_3]^{\zp}$ were given by Dickson \cite{dickson}.
Generators for  $\field[V_4]^{\zp}$ and $\field[V_5]^{\zp}$ can be found in \cite{shank}.
The problem of finding an explicit generating set for $\field[V_n]^{\zp}$ for $n>5$, remains open.
Even when the invariants of the indecomposable summands are well understood,
it can be difficult to construct generating sets for decomposable representations.
Campbell \& Hughes, in \cite{campbell-hughes}, describe a generating set for $\field[kV_2]^{\zp}$, proving
a conjecture of Richman \cite{richman1}. Generating sets have recently been constructed for $\field[V_2+V_3]^{\zp}$ \cite{sw:noeth}
and $\field[2V_3]^{\zp}$ \cite{cfw}.
Using the known generating sets, it is easy to see that $\beta(V_2)=\beta(V_3)=\beta(2V_2)=p$.
It follows from \cite{campbell-hughes} and \cite{richman1}, that 
$\beta(kV_2)=k(p-1)$ for $k>2$. By \cite[Remark~3.3]{se-sh}, $\beta(V_4)=2p-3$. This paper is devoted to
computing the Noether number for all remaining modular representations of $\zp$. 
An $\field\zp$-module is called {\it reduced} if it is has no trivial summands.
Note that, for any $\field\zp$-module $W$, we have $\beta(kV_1+W)=\beta(W)$, i.e., adding trivial summands does not change the Noether number.
Therefore, it is sufficient to compute $\beta$ for reduced $\field\zp$-modules. 
Also note that $V_n^{\zp}$, the vector space of fixed points, has dimension one.
Therefore, for any $\field\zp$-module $W$, $\dim(W^{\zp})$ is the number of indecomposable 
summands in the decomposition of $W$. 

\begin{proposition}\label{main_prop}
Suppose that $W$ is a reduced finite dimensional $\field\zp$-module.\\
(a) If $W$ contains a summand isomorphic to $V_n$ with $n>3$, then
$$\beta(W)=(p-1)\dim(W^{\zp})+p-2.$$ 
(b) If $W$ is isomorphic to $mV_2+\ell V_3$
with $\ell>0$, then $$\beta(W)=(p-1)\dim(W^{\zp})+1.$$
\end{proposition}
\begin{proof}
The first result follows from the lower bound given in Theorem~\ref{lbthm} and the upper bound given in Corollary~\ref{ubcor}.
The second result follows from the lower bound given in Theorem~\ref{lbthm} and the upper bound given in Corollary~2.8 of \cite{se-sh}.
\end{proof}

We note that the ``$2p-3$ conjecture''\cite[Conjecture~6.1]{sw:noeth} is a special case of Proposition~\ref{main_prop}~(a),
and that Conjecture~1.1 of \cite{se-sh} follows from Theorem~\ref{upper bound theorem}.

For general facts concerning the invariant theory of finite groups, see 
\cite{benson}, \cite{de-ke} or \cite{ne-sm}.

\section{Lower bounds} \label{lbsec}

It is well known that $\beta(V_2)=\beta(V_3)=\beta(2V_2)=p$. For all other reduced representations, $W$, the number
$(p-1)\dim(W^{\zp})$ is a lower bound for $\beta(W)$ (see \cite[Proposition~3.1]{sw:noeth}).
If $W=kV_2$ with $k>2$, then this lower bound is sharp (see \cite{campbell-hughes} and \cite{richman1}).
This section is devoted to proving the following theorem.

\begin{theorem} \label{lbthm} Suppose that $W$ is a  non-trivial reduced $\field\zp$-module. Define $k:=\dim(W^{\zp})$.
If  $W$ contains a summand isomorphic to $V_n$ with $n>3$, then $\beta(W)\geq k(p-1)+p-2$.
If  $W$ contains a summand isomorphic to $V_n$ with $n>2$, then $\beta(W)\geq k(p-1)+1$.
\end{theorem}

Recall that for a submodule $U\leq W$, $\beta(U)\leq\beta(W)$ \cite[Theorem~4.2]{sw:noeth}.
Thus to construct a lower bound, it is sufficient to find a lower bound for a submodule.
To prove the first result, we give a lower bound for $\beta((k-1)V_2+V_4)$ and to prove the second result,
we give a lower bound for $\beta\left((k-1)V_2+V_3\right)$. 

Choose a basis $\{x_1,y_1,x_2,y_2,\ldots,x_{k-1},y_{k-1},w,x,y,z\}$ for $\left((k-1)V_2+V_4\right)^*$
with $\Delta(y_i)=x_i$, $\Delta(x_i)=0$, $\Delta(z)=y$, $\Delta(y)=x$, $\Delta(x)=w$ and $\Delta(w)=0$.

\begin{lemma}\label{V4_lower} $\tr\left((y_1\cdots y_{k-1}z)^{p-1}y^{p-2}\right)$ is indecomposable in $\field[(k-1)V_2+V_4]^{\zp}$.
\end{lemma}
\begin{proof}
Define $H:=\tr\left((y_1\cdots y_{k-1}z)^{p-1}y^{p-2}\right)$.
The proof is by induction on $k$. For $k=1$, the result follows from Remark~3.3 of \cite{se-sh}.
Assume $k>1$. Suppose, by way of contradiction, that $H=\sum f_ih_i$ for some
$f_i,h_i\in\field[(k-1)V_2+V_4]_+^{\zp}$.
We may assume that $f_i$, $h_i$ and $f_ih_i$ are homogeneous with respect to multidegree,
and that each $f_ih_i$ has multidegree $(p-1,p-1,\ldots,p-1,2p-3)$.
Use the inclusion of $V_1+(k-2)V_2+V_4$ into $(k-1)V_2+V_4$ to define a projection
$\pi:\field[(k-1)V_2+V_4]\to\field[V_1+(k-2)V_2+V_4]$. Note that $\pi$ is the algebra homomorphism
taking $x_1$ to $0$, $y_1$ to $x_1$ and fixing the other variables. 
Since $\pi$ is equivariant, 
$\pi(H)=\tr\left(\pi\left((y_1\cdots y_{k-1}z)^{p-1}y^{p-2}\right)\right)=x_1^{p-1}\tr\left((y_2\cdots y_{k-1}z)^{p-1}y^{p-2}\right)$.
Collecting all factors of $x_1$, write $\pi(f_i)=x_1^{n_i}\widetilde{f_i}$ and $\pi(h_i)=x_1^{m_i}\widetilde{h_i}$.
Using the homogeneous multidegree assumption, $n_i+m_i=p-1$ and 
$\widetilde{f_i}$,$\widetilde{h_i}\in\field[(k-2)V_2+V_4]^{\zp}$. 
Thus $\pi(H)=x_1^{p-1}\sum \widetilde{f_i}\widetilde{h_i}$.
Therefore $\widetilde{H}:=\tr\left((y_2\cdots y_{k-1}z)^{p-1}y^{p-2}\right)=\sum \widetilde{f_i}\widetilde{h_i}$.
By the induction hypothesis, $\widetilde{H}$ is indecomposable in $\field[(k-2)V_2+V_4]^{\zp}$.
It is clear that $ \widetilde{f_i},\widetilde{h_i}\in \field[(k-2)V_2+V_4]^{\zp}$. 
Suppose that one of the factors, say $\widetilde{f_1}$, is of degree zero. Then the multidegree of $f_1$ is
$(n_1,0,0,\ldots,0)$ with $n_1<p$. Hence $f_1$ is a homogeneous element of degree less than $p$ in 
$\field[x_1,y_1^p-x_1^{p-1}y_1]\cong\field[V_2]^{\zp}$. Thus $f_1=cx_1^{n_1}$ for some $c\in\field$.
Therefore $\pi(f_1)=0$, giving $\widetilde{f_1}=0$. Hence $\sum \widetilde{f_i}\widetilde{h_i}$
determines a decomposition of $\widetilde{H}$, providing the required contradiction.
\end{proof}

Choose a basis $\{x_1,y_1,\ldots,x_{k-1},y_{k-1},x,y,z\}$ for $\left((k-1)V_2+V_3\right)^*$
with $\Delta(y_i)=x_i$, $\Delta(x_i)=0$, $\Delta(z)=y$, $\Delta(y)=x$ and $\Delta(x)=0$.
We use the graded reverse lexicographic order with $x_1<y_1<x_2<\cdots <y_{k-1}<x<y<z$. 

\begin{lemma}
$\tr\left((y_1\cdots y_{k-1}z)^{p-1}y\right)$  is indecomposable in $\field[(k-1)V_2+V_3]^{\zp}$.
\end{lemma}

\begin{proof}
The proof is by induction on $k$. For $k=1$, $\field[V_3]^{\zp}$ is a hypersurface
with generators in degrees $1$, $2$, $p$ and $p$ (see \cite{dickson} or \cite[Proposition~5.8]{peskin}); the element $\Tr(yz^{p-1})$ has lead monomial
$y^p$ and may be chosen to be one of the generators. The induction step is essentially the same as the induction step in the
proof of Lemma~\ref{V4_lower}.
\end{proof}

\section{An upper bound}\label{ubsec}

Suppose that $W$ is a reduced finite dimensional $\field\zp$-module. The purpose of this section is to
construct an upper bound on $\beta(W)$. Decompose $W$ into a direct sum of indecomposable modules:
$W=\oplus_{i=1}^k W_i.$ For each $i$, choose an $\field\zp$-module generator $z_i$ for $W_i^*$. Define $d_i:=\dim(W_i)$.
Then $$\bigcup_{i=1}^k\{\Delta^t(z_i)\mid t=0,1,\ldots,d_i-1\}$$ is a basis for $W^*$. We use this basis to determine 
generating variables for $\field[W]$. Let $A$ be the subalgebra of $\field[W]$ generated by
$$\bigcup_{i=1}^k\{\Delta^t(z_i)\mid t=1,\ldots,d_i-1\},$$
i.e., the subalgebra formed by omitting the ``terminal variables''.
Suppose that $m=u_1u_2\cdots u_{p-1}$ is a monomial of degree $p-1$ in $A$. For each $u_i$, there exists a variable, say $w_i$,
with $\Delta(w_i)=u_i$. Define $m'=w_1w_2\cdots w_{p-1}$. Note that $m'$ is a monomial of degree $p-1$ in $\field[W]$.
Define
$$F:=\sum_{\ell\in \fp}\prod_{j=1}^{p-1}\left(w_j-\sigma^{\ell}(w_j)\right),$$
where $\fp$ denotes the field with $p$ elements.
For $S\subseteq\{1,2,\ldots,p-1\}$, define $S'=\{1,2,\ldots,p-1\}\setminus S$. Further define 
$X_S=\prod_{j\in S}w_j$ and $X_{S'}=\prod_{j\in S'}w_j$.

\begin{lemma}
$$F=\sum_{S\subseteq \{1,2,\ldots, p-1\}}(-1)^{|S|}X_{S'}\tr(X_S).$$
\end{lemma}

\begin{proof} Expanding, and using the fact that $\sigma^{\ell}$ is an algebra automorphism, gives
$$\prod_{j=1}^{p-1}\left(w_j-\sigma^{\ell}(w_j)\right)=\sum_{S\subseteq \{1,2,\ldots, p-1\}}(-1)^{|S|}X_{S'}\sigma^{\ell}(X_S).$$
Summing over $\ell\in \fp$ gives the required formula.
\end{proof}

We use a graded reverse lexicographic order with $\Delta^t(z_i)>\Delta^{t+1}(z_i)$ for $t<d_i-1$, and denote the lead term of a polynomial $f$ by $\lt(f)$.

\begin{lemma}\label{mufit} $\lt(F)=-m$.
\end{lemma}
\begin{proof}
Note that $\sigma^{\ell}(w_j)=w_j+\ell u_j +\cdots$, where the missing terms are lower in the order.
Thus the lead term of $w_j-\sigma^{\ell}(w_j)$ is $-\ell u_j$.
Since the lead term of a product is the product of the lead terms, we have
$\lt\left(\prod_{j=1}^{p-1}\left(w_j-\sigma^{\ell}(w_j)\right)\right)=(-\ell)^{p-1}\prod_{j=1}^{p-1}u_j=\ell^{p-1}m$.
Thus $F$ is a sum of polynomials  all with lead monomial $m$. The result follows from the fact that $\sum_{\ell\in\fp}\ell^{p-1}=-1$.
\end{proof}

Recall that the {\it Hilbert ideal} of $W$ is the ideal in $\field[W]$ generated by homogeneous invariants of positive degree.
The ring of {\it coinvariants}, denoted by $\field[W]_{\zp}$, is the finite dimensional graded algebra consisting of the
quotient of $\field[W]$ by the Hilbert ideal.

\begin{theorem}\label{upper bound theorem} The top degree of $\field[W]_{\zp}$ is bounded above by $k(p-1)+p-2$.
\end{theorem}
\begin{proof}
The Hilbert series of a graded ideal in $\field[W]$, and the series for the corresponding ideal of lead terms, coincide. 
Thus it is sufficient to bound
the top degree of the quotient by the lead term ideal of the Hilbert ideal.
It follows from Lemma~\ref{mufit} that every monomial of degree $p-1$ in $A$ is in the
lead term ideal. Furthermore, $z_i^p$ is the lead term of the orbit product of $z_i$. Thus the monomials lying outside the lead
term ideal are of the form $\gamma\prod_{i=1}^kz_i^{n_i}$ were $\gamma$ is a monomial in $A$ with degree at most $p-2$, and $n_i\leq p-1$.
\end{proof}

\begin{corollary}\label{ubcor}
The Noether number of $W$ is bounded above by $k(p-1)+p-2$.
\end{corollary}
\begin{proof}
The top degree of the coinvariants gives an upper bound for the degrees of a generating set of the image of transfer.
It follows from \cite[Lemma~2.12]{hk} that $\field[W]^{\zp}$ is generated by transfers, orbit products and elements
of at most degree $k(p-1)-(\dim(W)-k)$. Thus the upper bound for the top degree of the coinvariants
given in Theorem~\ref{upper bound theorem} gives an upper bound for the Noether number of $W$.
\end{proof}

\end{document}